\documentclass[12pt,a4paper,reqno]{amsproc}
\usepackage[colorlinks]{hyperref}

\usepackage{amsthm,array,mathrsfs}
\usepackage{amsmath,amssymb,enumitem}

\usepackage{tikz}
\usetikzlibrary{patterns}
\usetikzlibrary{calc}

\newtheorem{theorem}{\rm\bf Theorem}[section]

\newtheorem*{theorem*}{Theorem}
\newtheorem*{theorem 1}{\rm\bf Proposition 1}
\newtheorem*{theorem 2}{\rm\bf Proposition 2}
\newtheorem*{conj*}{Conjecture}

\theoremstyle{definition}

\theoremstyle{remark}
\newtheorem{remark}[theorem]{\rm\bf Remark}

\def\half#1#2{\begin{matrix}\frac{#1}{#2}\end{matrix}}

\def\R#1{\mathbb{R}^{#1}}

\def\scal#1#2{\langle #1; #2 \rangle}

 \DeclareMathOperator{\re}{Re}
\DeclareMathOperator{\im}{Im}

\DeclareMathOperator{\trace}{tr}

\begin{document}

\title{V.M. Miklyukov: from dimension $8$ to nonassociative algebras}

\author{Vladimir G. Tkachev}
\address{Department of Mathematics, Link\"oping University, Link\"oping, 58183, Sweden}
\email{vladimir.tkatjev@liu.se}

\begin{abstract}
In this short survey we give a background and explain some  recent developments in algebraic minimal cones and  nonassociative algebras. A good deal of this paper is recollections of my collaboration with my teacher, PhD supervisor and a colleague, Vladimir Miklyukov on minimal surface theory that motivated the present research. This paper is dedicated to his memory.
\end{abstract}
\dedicatory{In memory of V.M.~Miklyukov (1944-2013)}
\maketitle




\section{Introduction}
In this short note, I try to explain  how and why certain nonassociative algebra structures arise in the context of minimal cones. The starting point of this relatively recent subject comes from the Bernstein property for minimal graphs. More precisely, it comes from an attempt to understand the breakdown in higher dimensions $n\ge8$ of the celebrated S.N.~Bernstien theorem from 1915 on minimal two-dimensional graphs. Although the original question remains still unanswered, it is clear now that the nonassociative aspect of certain classes of elliptic type PDEs including eiconal equation, minimal surface equation and even more general classes (see \cite{NTVbook}, \cite{Fox19a}) is relevant and not a coincidence. Remarkably, very related classes of nonassociative algebras appear very naturally in a very different,  group theoretic context, see \cite{HRS15}, \cite{Rehren17}, \cite{Ivanov09} and the references thererin.

I do not set myself the goal to cover all the latest achievements in this area or to give a  complete overview. The subject is still under development and many basic questions are still open.  The interested reader is referred to the recent monograph \cite{NTVbook} for unexpected  connections of nonassociative algebras to regularity theory of fully nonlinear PDEs. See also the recent papers \cite{Tk14}, \cite{Tk18a}, \cite{Tk18b}, \cite{KrTk18a}, \cite{Tk18e}, \cite{Tk19a} for further results and basic concepts considered below.

\textbf{Acknowledgments.} I thank Prof. A.V.~Loboda for his careful reading of the manuscript and helpful suggestions.

\section{Bernstein's problem}
Let me first explain some relevant historical details.
Miklyukov has been my supervisor ever since I was a second year student, leading and coaching me all the way through my PhD. I learned so much from him, both academically and professionally. Miklyukov always had a very extensive range of interests covering several domains in nonlinear and geometrical analysis. Somewhere in the very beginning, he once told me that around 1975-1976, one of his colleagues suggested that he deal with minimal surfaces and try to apply there the previously developed methods and results from the theory of quasiconformal mappings. The idea turned out to be extremely fruitful and already in 1977-80 he published a series of important papers on boundary and asymptotic behaviour, Liouville type theorems for a large class of quasilinear euqations \cite{Mik77a}, \cite{Mik78a}, \cite{Mik79a}, \cite{Mik79b}, \cite{Mik80a}. In these papers, he masterly developed several concepts and results from the classical complex analysis, potential theory and elliptic type PDEs, and   applied the new methods, in particular,  to study of minimal surfaces in the Euclidean space. In summary, his approach  was based on a virtuous interplay of several  key ingredients including the extremal length, variational capacity, fundamental frequency, uniformization theory, Beltrami equation, Liuville, Fragm\'{e}n-Lindel\"of,  Denjoy-Carleman-Ahlfors and Wiman type theorems. One of the main goals of the proposed by Miklyukov approach was the famous Bernstein theorem and its generalizations.

Since an important role of this result playing in our further context, I briefly recall some relevant definitions and concepts. A surface in $\R{3}$ is called minimal if its mean curvature vanishes everywhere. From an analytic point of view, this is equivalent to saying that a minimal surface is a stationary point of the area functional. This implies that locally any minimal surface satisfies a very nice quasilinear equation of the second order. In 1903, Sergei Natanovich Bernstein, published a short note in \textit{Comptes Rendus} containing some crucial results on the analyticity solutions of  second order elliptic partial differential equations, hereby solving (in the $C^3$-regularity class) the 19th problem addressed by Hilbert at the First International Mathematical Congress in 1900. Between 1903 and 1918, S.N.~Bernstein published several important memoirs on the regularity and a priori estimates for quasilinear elliptic PDEs with further applications to surfaces with a prescribed mean curvature. One of the most remarkable results obtained by Bernstein was his famous theorem of 1915 (first published in \textit{Communications of the Kharkov Mathematical Society}) asserting that any \textit{entire}, i.e. defined in the whole plane $\R{2}$, solution $u(x)$ of the minimal surface equation
\begin{equation}\label{MS}
\mathrm{div}\, \frac{Du}{\sqrt{1+|Du|^2}}=0, \qquad x\in \R{2},
\end{equation}
must be an affine function, i.e. $u=ax+by+c$.

The Bernstein theorem is remarkable in many aspects.
Unlike the classical Liouville result on bounded harmonic functions, it claims that a solution must be `trivial' (affine) without any additional assumptions on the growth of a solution. Some possible explanation of the latter phenomenon follows from the fact that the minimal surface equation is much more symmetric than the Laplace equation. Indeed, equation \eqref{MS} is invariant under the full orthogonal group $O(3)$ (rotations of $\R{3}$), while the property of being a harmonic function is invariant under the action of a smaller subgroup $O(2)$ of rations of $\R{2}$, i.e. a rotation of a graph of a harmonic function in $\R{3}$ is no longer harmonic. In particular, the large symmetry group of \eqref{MS} makes it possible to associate  to any solution of \eqref{MS} another \textit{a priori} bounded solution. More precisely, Bernstein remarks that  $v(x)=\arctan \frac{\partial f}{\partial x_1}$ is always a (bounded!) solution to an elliptic linear equation. Then the claim follows by application of the Liouville theorem proved by Bernstein four years earlier.

Unfortunately, this method does not give clues about possible generalizations of the Bernstein result onto $n$-dimensional minimal hypersurfaces in $\R{n+1}$. The story of the Bernstein theorem for $n\ge 2$ and searching of counterexamples requires completely different approaches and it  is a fascinating part of the mathematical history including several brilliant names as W.H.~Fleming, H. ~Federer, E. De Giorgi, F.J.~Almgren, M.~Miranda, J. Simons, E.~Bombieri, E.~Giusti, R.~Osserman, J.C.C.~Nitsche, L.~Simon. We refer the interested reader to a recent surveys of Leon Simon \cite{SimonV} and Mario Miranda \cite{Miranda} for the further reading and more details. The result summarizing several papers states that the \textit{Bernstein property} (i.e. the claim that an entire solution to a certain second order elliptic type PDE is an affine function) holds true for all dimensions $2\le n\le 7$ and in the higher dimensions $n\ge8$, there are nontrivial solutions of \eqref{MS} over the whole $\R{n}$.

Although Bernstein’s problem has already been settled 50 years ago, it remains an important cornerstone of analysis and geometry, a kind of incomprehensible and unattainable beautiful Everest with immortal charm, attracting the constant attention of specialists. And although we have today many answers, the key question remains unanswered: why does the Bernstein property collapse in higher dimensions?

In this respect, the two-dimensional case takes a very special place. It is intimately related to the fact that the two-dimensional Euclidean space has a natural complexi\-fi\-cation:
\begin{equation}\label{Complex}
\R{2}\cong \mathbb{C}^1.
\end{equation}
The existence of the complex structure on $\R{2}$ implies variety of different versions of the proof and approaches to the Bernstein theorem in two dimensions. On the other hand, the complex structure also implies the existence of the Weierstrass-Enneper parameterization and, as a corollary, an exceptional variety of species  in  the two-dimensional minimal surface zoo. More precisely, it has been proved by  Weierstrass and Enneper (around 1863) that locally any two-dimensional minimal surface is given by the real part of a certain holomorphic curve. In other words, it can be parametrized by
\begin{equation}\label{wier}
x(z)=(\re \phi_1(z), \re \phi_2(z), \re \phi_3(z)),
\end{equation}
where the triple of meromorphic functions $\phi_i(z)$ on a Riemann surface satisfies
$$
\sum_{i=1}^3\phi_i'^2(z)=0.
$$
Conversely, any such a triple generates a minimal surface. This representation makes the similarity between harmonic (and meromorphic) functions in $\R{2}$ and minimal surfaces in $\R{3}$ more rigorous. For example, setting $\phi_3=0$ yields $\phi_2=\sqrt{-1}\phi_1(z)$, hence
$$
x(z)=(\re \phi_1(z), -\im \phi_1(z),0)
$$
becomes essentially the Cauchy-Riemann representation expressing that  any harmonic function $u(x_1,x_2)$ in a plane can be locally written as the real part of a holomorphic function.

\section{Nevanlinna theory on minimal surfaces}
In the fall 1982,  being a 2nd-year student, I started a project with Miklyukov on properly immersed minimal surfaces in $\R{3}$. Recall that the Bernstein theorem requires that a minimal surface must be a graph over a whole plane. Miklyukov's idea was to relax the latter condition and replace it by an appropriate projection property, in other words, to find a `quantitative' Bernstein property. Miklyukov has always had incredible analytical and geometrical intuition based, in particular,  on his previous research on the regularity of quasiconformal mappings and the classical function theory.

As an appropriate analogue of the projection property, Miklyukov thought to develop the concepts of the counting function and the defect relations in the classical Nevanlinna theory of meromorphic functions to immersed minimal surfaces. Recall that the Nevanlinna theory describes the asymptotic distribution of solutions of the equation $f(z) = a$, $f(z)$ being a meromorphic function,  as the value $a$ varies. Then the counting function is just the logarithmic average number of solutions. A fundamental tool of the theory is the Nevanlinna characteristic which accumulates both the counting function with the growth function (or the proximity function). Taking into account the made above remark on the similarity between holomorphic functions and minimal surfaces, the idea seemed very natural.

The counting function in the  context of minimal surfaces is just the normalized average over a circle of radius $r>0$ of the multiplicity of the orthogonal projection of a minimal surface onto a fixed plane in $\R{3}$. Thus, for the graph over a plane, the counting function is just a  constant equal to 1. Although the very close similarity between meromophic functions and minimal surfaces, it was absolutely unclear how to adopt the analytical arguments for functions (which have an algebra structure) to surfaces. A promising approach was to try mixup the harmonicity of the coordinate functions on a minimal surface with some standard tools of the classical function theory: the extremal length, variational capacity, the co-area formula\footnote{Miklyukov used to refer to the coarea-formula as the \textit{Kronrod-Federer formula}, according to the original Russian edition of the book of Burago and Zalgaller \cite{Burago}.} and the fundamental frequency.

The capacity and fundamental frequency are quite standard instruments in potential analysis and PDEs, while the concept of the extremal length comes from  Teichm\"uller theory of quasiconformal mappings and it is not as standard  in PDEs. These tools has been recently employed by Miklyukov  in his elegant approach  \cite{Mik79a}  to the Bernstein problem in two dimensions. Some preliminary results obtained by him in this and a couple of others papers suggested in particular that it is very  natural to expect a `relaxed'  Bernstein property for minimal quasigraphs. The idea was to connect the surface geometry to the counting function by the well-known identity between extremal length and the conformal capacity, since the latter could be effectively estimated on surfaces.

The first part of this project has been finished and submitted in May 1984, and published three years later in \cite{MT87}\footnote{Submission and publication of this paper coincided with the most difficult period of Miklyukov, between 1985 and 1989, who lost his position as Chair of Mathematical Analysis and Function Theory (the chair itself has been abolished in 1985), it is a long different story. Miklyukov was returned to the title of professor in 1988, and the chair was returned its name only in 1992. It took almost four years to publish our paper \cite{MT87}, there were two negative reviews written by `colleagues', but then the situation changed mysteriously and we got a positive review. During that period, I finished the last year of university, then served a year and a half in the army, and even entered my PhD program.}, we were able to prove two results for different types of projections. Here is the \href{https://mathscinet.ams.org/mathscinet/search/publdoc.html?arg3=&co4=AND&co5=AND&co6=AND&co7=AND&dr=all&pg4=AUCN&pg5=TI&pg6=MR&pg7=ALLF&pg8=ET&r=1&review_format=html&s4=&s5=&s6=MR0912614&s7=&s8=All&sort=Newest&vfpref=html&yearRangeFirst=&yearRangeSecond=&yrop=eq}{review} from MathSciNet:

\begin{footnotesize}
\begin{quote}
This article concerns the absence of nontrivial noncompact parametric minimal surfaces in the Euclidean $3$-space satisfying some additional geometrical properties. A typical known result was the following. Let $x:M\to \R{3}$ be a $C^2$ minimal imbedding of a noncompact orientable $2$-manifold. If $x(M)$ lies between two parallel planes and if the induced metric on M is complete, then $x(M)$ is necessarily a plane ([L. P. de M. Jorge and F. V. Xavier, Ann. of Math. (2) 112 (1980), no. 1, 203–206; \href{https://mathscinet.ams.org/mathscinet/search/publdoc.html?r=1&pg1=MR&s1=584079&loc=fromrevtext}{MR0584079}).  Does the conclusion still hold provided $x(M)$ only lies in a half-space (a problem stated by Calabi)? In this connection, the authors assume that $x$ is proper and $x(M)$ lies in a half-space, say in $\{x_3>0\}$. Properness implies that $x(M)$ is topologically complete (in $\R{3}$), which is weaker than metric completeness. They also assume that the average number of intersection points with $x(M)$ on a line parallel to the $x_3$-axis passing through a generic point in $\{x_3=0\}$ situated at distance t from the origin is $o(\log t)$. They deduce the triviality of $x$, via a "parabolicity'' property of $M$ (equipped with the metric induced by $x$).

   \begin{flushright}
   Reviewed by Philippe Delano\"e
   \end{flushright}
   \end{quote}
   \end{footnotesize}

\section{The Blue Notebook and the dimension $8$}
Another important motivation for that project was that Miklyukov aimed to achieve a more conceptual explanation of the Bernstein property and to find a way to approach minimal submanifolds in higher dimensions. The project required me to read extensively from various fields including integral geometry, elementary geometric measure theory, value distribution theory of meromorphic functions,  manifold theory  and Grasmannian geometry. At that time  I met with Miklyukov's `blue notebook'.

It was an old-fashioned thick A5 notebook, with  embossed \textsc{Tpect} on the cover and filled with  translations of various articles from English and German. One of the first was a translation of a paper of Federer on geometric measure theory.  Miklyukov showed me a translation made in a red ink of the famous Milnor's \emph{On manifolds homeomorphic to the 7-sphere} \cite{Milnor56}. The seven-page Milnor amazing proof was a virtuous combination of the Morse theory,  Thom's cobordisms, and the  quaternionions. He constructs his famous exotic 7-dimensional  spheres  $M_k^7$ as $S^3$-bundles over $S^4$ in a very explicit  way using quaternionic multiplication. The Morse theory then  applies that the constructed $M_k^7$ must be homeomorphic to the standard sphere $S^7$, while the $\lambda$-invariant may be different for distinct values of $k$. The latter implies the existence of different differentiable structures on $S^7$.

I remembered this moment, because it was rather unusual situation. Miklyukov told me that he believe that the breakdown of the Bernstein property in dimension 8 could  have a connection with the existence of the exotic structures for higher dimensional spheres and the Bott periodicity. He sketched his motivation based on both formal and intuitive links and subtle analogies between these two completely different worlds\footnote{Much latter I found an interesting collection of around 80 different sources, papers and links on `Eight in algebra, topology and mathematical physics' on Andrew Ranicki’s homepage \cite{Ranicki}.}. Some promising bridge could be the $N$-averages of the fundamental frequency defined and studied in \cite{Mik80a}:
\begin{equation}\label{lamm}
\lambda_{\alpha}(D, N)=\inf\{\frac{1}{N}\sum_{i=1}^N\lambda_{\alpha}(D_i):\,\sqcup_{i=1}^N D_i \subset D\},
\end{equation}
where $D$ is an open domain  of a Riemannian manifold $M^n$. One can easy to show that $\lambda_{\alpha}(D, N)$ is increasing in $N$. For $M^n=\R{n}$ the structure of $\lambda_{\alpha}(D, N)$ is well understood and an optimal lower estimate
$$\lambda_{\alpha}(D, N)\ge c_n(N/|D|)^{1/n}
$$
holds. The most interesting nontrivial case of the Euclidean spheres $S^{n-1}\subset \R{n}$ is more subtle. There are still some analogues of the above lower estimate but they are not optimal with respect to $N$. The question is already very nontrivial for the two-dimensional sphere.\footnote{I just remark that another very context where $\lambda_{\alpha}(D, N)$ naturally appear is to the nodal level sets of the $k$th eigenfunction of the Laplacian operator on compact manifolds. However, a connection between $k$ and $N$ is unclear for general $n\ge2$.}

The essence of $\lambda_{\alpha}(D, N)$ comes from the following crucial observation made in \cite{Mik80a}: a sharp lower estimate for the average fundamental frequency of spheres implies sharp versions of Liouville, Denjoy-Carleman-Ahlfors and Wiman type theorems for higher-dimensional subsolutions of quasilinear ellitptic type equations on the corresponding manifold. The minimal surface equation fits naturally  in this class. The link between the existence of exotic spheres is more subtle but Miklyukov believed that it could be related to the \textit{extremal} combinatorial structure for the variational problem \eqref{lamm}.

A noteworthy thing was that the idea had a clear algebraic tone, a rather unusual shade for Miklyukov  who seemed  primarily an analyst and a geometer. It was quite obvious that this idea was absolutely important and close for him. But the problem seemed much too
difficult and was actually a pure rhetoric problem: we had no idea how to get started.


\section{Cubic cones enter the picture}
Between 2003 and 2010, while a researcher at Royal Institute of Technology (Stockholm) I worked together B\"orn Gustafsson on a completely different project in complex analysis: the Hele-Shaw problem, complex moments and meromorphic resultants, see his nice review  \cite{Gust2014} about some of these developments. During Fall semester 2008, I switched back to minimal surfaces and began my work on the Hsiang problem. Some preliminary results were obtained during my train journeys between Stockholm and Uppsala, where I temporarily taught at the  Swedish National Graduate School in Mathematics and Computing (FMB).

My interest to the Hsiang problem has been  motivated by two principal questions: a better understanding of the algebraic nature of the known examples of entire graphs in dimensions $\ge 8$ and generalizations of an example of a fourfold periodic minimal hypersurfaces in $\R{4}$ I discovered at that time, see \cite{Tk13a}. This latter example was a logical sequel of the three-dimensional double-periodic examples constructed earlier together with Vla\-dimir Sergienko \cite{SergTk2000}, \cite{SergTk1999}. More precisely, the obtained fourfold periodic example was an embedded minimal hypersurface with a $D_4$-symmetry group and isolated singularities of the Clifford cone type. Until recently, the only known non-trivial (i.e. distinct from cones) examples of embedded minimal hypersurfaces in $\R{n+1}$ with \textit{finitely many} isolated singularities were constructed by Cafarelli-Hardt-Simon  \cite{CafHarSi}, N.~Smale \cite{SmaleN}, and Harvey-Lawson \cite{HL}. Thus, the above $D_4$-invariant minimal hypersurface was the first example of an embedded minimal hypersurface with \textit{infinitely many} isolated singular points. This naturally led to the following question: given a lattice  $L$ in $\R{n}$ and a fixed  minimal cone $K$, is it possible to construct an $L$-periodic embedded minimal hypersurface with isolated singularities of the type $K$ at the points of $L$?

One obvious obstacle to this task was a shortage of available examples of minimal cones. It is well-known that in any dimension $n\ge4$ there exist exactly $[\frac{n-2}{2}]$ non-congruent \textit{quadratic} minimal cones, all  classified in a landmark paper of Wu-yi~Hsiang \cite{Hsiang67} of 1967 published in the first issue of Journal of Differential~Geometry.  In the same paper, Hsiang  formulate several  problems on general real algebraic minimal cones and by using invariant theory constructs explicitly four new examples of cubic minimal cones in dimensions $n=5,8,9$ and $15$. All the obtained (and known so far irreducible) examples of cubic minimal cones satisfy the nonlinear 2nd order PDE
\begin{equation}\label{radial}
|D u|^2\Delta u-\half{1}{2}Du\cdot D|Du|^2=\theta |x|^2 u, \quad x\in \R{n},
\end{equation}
where $\theta\in \R{}$ is a structure constant. A cubic polynomial solution of  (\ref{radial}) is said to be a \textit{Hsiang eigencubic} \cite{Tk10c} (or \textit{radial eigencubic},   \textit{REC} for short, according \cite{Tk10c}).

In general, an algebraic minimal cone of degree $d$ is determined by a homogeneous polynomial solution $u$ of \eqref{radial} with a certain homogeneous polynomial in $x$ of degree $2d-4$ instead of the quadratic from $\theta |x|^2$ in the right hand side. The first nontrivial is the case of degree $d=3$ solutions of \eqref{radial}, but Hsiang remarks that `the algebraic difficulties involved in such a problem are rather formidable' and then asks to characte\-rize at least all cubic homogeneous polynomial solutions of \eqref{radial} \cite[p.~258, 265]{Hsiang67}.

During 2008--2010, I obtained a particular solution of the Hsiang problem by  applying a rather straight\-forward approach and summarized in the preprint \cite{Tk10b}. Among the results obtained there, I mention the following:
\begin{enumerate}[label=(\roman*)]
\item\label{item1} Any   Hsiang eigencubic $u$ is harmonic, unless it is trivial (i.e. congruent to the trivial one-variable cubic polynomial $ax_1^3$ under an isometry of $\R{n}$).
\item\label{item2} There is an infinite family of Hsiang eigencubics  (called Clifford type eigencubics) explicitly parameterized by symmetric Clifford systems.
\item\label{item3}
  There are only finitely many non-Clifford type eigencubics (the so-called \textit{exceptional} eigencubics), and the the list of a priori possible dimensions was established.
\item\label{item4} Any Hsiang eigencubic satisfies the cubic trace identity for its Hessian matrix: the trace of the third power of the Hessian  of $u$ is proportional to $u$, i.e.
    $$\trace (D^2u)^3=au,
    $$
    where $a\in \R{}$ contains a certain intrinsic invariant information about $u$.
\item\label{item5} An eigencubic $u$ is exceptional if and only if it satisfies the second order identity:
    $$
    \trace (D^2u)^2=c|x|^2.
    $$
    In fact, there are four eigencubics (coming from trialities) that are intermediate between Clifford type and exceptional, the so-called mutants, which satisfy both the quadratic trace identity and Clifford type representation.
\item\label{item6} A new example of an exceptional eigencubic in dimension $21$ was constructed by using the octonion algebra. The corresponding Hsiang eigencubic is simply the following real part:
    $$
    u(x)=\re (w_1w_2w_3),
    $$
    where $x=(w_1,w_2,w_3)$ with $w_i$ being three independent imaginary octonions in $\im \mathbb{O}\cong \R{7}$.
\end{enumerate}

Some remarks are in order here. First, the harmonicity of  nontrivial Hsiang eigen\-cubics is a rather striking property which means that any eigencubic  satisfies in fact to a system of two  second order PDEs. Unfortunately, the available proofs of this property do not shed light on a conceptual understanding why this extra PDE constraint does actually hold.  It is   unclear neither if an analogous phenomenon holds for eigencubics of higher degrees. Next, the definition of an exceptional Hsiang eigencubic is somewhat negative, thus  nonconstructive, hence it was very desirable to find any constructive way to distinguish exceptional eigencubics from eigencubics of  Clifford type. In this respect, the result \ref{item5} plays the fundamental role.

The finiteness result \ref{item3} is very striking and one immediately recognizes here a parallel between the above dichotomy (Clifford vs exceptional) and classical  dichotomies `regular vs exceptional (or sporadic)' in simple Lie algebras, finite simple groups, ADE-classification of singularities etc. Even we know today that Hsiang eigencubics have close connections with Clifford and Jordan algebras, and, thus, to classification of simple Lie groups, this phenomenon is not completely clear and requires a further study. I mention  also an interesting interplay between \ref{item5} and the Killing Einstein property studied very recently in important works of Daniel J.F. Fox \cite{Fox19a}.

Finally note that both \ref{item4} and \ref{item5}  provide us with simple effective  algorithmic tools for identifying an eigencubic and its  type (in particular, the constant $c$ in \ref{item4} contains a lot of intrinsic information about $u$ and its geometry).

The above results were partially published in \cite{Tk10c} (the classification of  the Clifford type examples in \ref{item2}) and  in \cite{Tk10a}  (a preparatory key result for the proof of \ref{item3} generalizing a Cartan theorem has been established).

The proof of the finiteness result \ref{item3}  relies on the famous Hurwitz-Radon theorem which characterizes the possible dimensions of composition of quadratic forms \cite{Shapiro}. More precisely, the classical  Euler four-square and Degen eight-square identities generalize the well-known multiplicative relationships between sums of squares in two variables
$$
(x_1^2+x_2^2)(y_1^2+y_2^2)=(x_1y_1-x_2y_2)^2+(x_1y_2+x_2y_1)^2.
$$
Then the Hurwitz-Radon theorem asserts that an identity of the kind
$$
(x_1^2+\ldots +x_p^2)(y_1^2+\ldots +y_q^2)=b_1(x,y)^2+\ldots +b_p(x,y)^2
$$
is possible for some bilinear forms $b_i$ if and only if $q=\rho(p)$, where  the Hurwitz-Radon function $\rho$ is defined by
\begin{equation}\label{foll}
\rho(m)=8a+2^b, \qquad \text{if} \;\,\,m=2^{4a+b}\cdot \mathrm{odd} , \;\; 0\leq b\le 3.
\end{equation}
This in particular implies the finiteness of exceptional eigencubics and also imposes obstructions  on their possible dimensions.

\begin{remark}
 In November 2010, I learned from Prof.~Zizhou Tang on the 1993 paper\footnote{I would like to thank Prof.~Yan Wenjiao for sending me a copy of this paper.} (in Chinese) \cite{Peng} of Peng Chin-Kuei and Xiao Liang where they proved \ref{item2}-\ref{item3}  using a similar approach but under an additional assumption that an eigencubic is harmonic, i.e. a priori requiring \ref{item1} as a condition. To my  knowledge this is the only work on the Hsiang problem which provides a particular classification similar to the above.
\end{remark}

\section{Nonassociative algebras of cubic forms}
Thus, the most difficult part of the Hsiang problem is  to determine \textit{which of the feasible dimensions are actually realizable}. For anyone who tried to deal with this problem, it became clear that at a certain stage any further progress faced numerous analytical and algebraic obstacles, thus an alternative, more transparent and conceptual, approach was required. To achieve a complete classification and better understanding of exceptional eigencubics, the straightforward approach of \cite{Tk10b} was insufficient. It gives no good idea why the Clifford and Jordan algebra structures may arise in the context of cubic minimal cones, or at least in the context of the Hsiang equation \eqref{radial}.

On the other hand, some hints were already in the above results. First, an easy calculation reveals  that any cubic form satisfying the Cartan-M\"unzner system \begin{equation}\label{Muntzer0}
|D u(x)|^2=9|x|^{4},
\end{equation}
\begin{equation}\label{Muntzer1}
\Delta u(x)=0
\end{equation}
 also satisfies the radial eigencubic equation \eqref{radial}. Recall that the latter system appears naturally in the context of isoparametric\footnote{A hypersurface is called isoparametric if it has constant principal curvatures \cite{Cecil1}} hypersurfaces in the spheres.  Cartan in 1938 established in \cite{Cartan38} that there exists exactly four congruency classes of such solutions.

 These \textit{Cartan isoparametric cubics} play a dual role in the classification \ref{item1}--\ref{item6}. First, they are simplest exceptional eigencubics corresponding to the Peirce parameter $n_2=0$. On the other hand,  the restriction of \textit{any} eigencubic with $n_2\ne0$ on a certain subspace is always a Cartan isoparametric cubic.
The latter property together with the classification result in \cite{Tk10a} were key ingredients in the proof of the finiteness result \ref{item3} above.

Another remarkable observation is that any Cartan isoparametric cubic is the \textit{generic norm} in the trace-free subspace of an appropriate simple rank 3 formally real Jordan algebra (see for example  \cite{BS2}). The simplest model of a Jordan algebra (and actually one half of possible examples) is the algebra of all real $n\times n$ symmetric matrices with the Jordan multiplication
$$
A\circ B=\frac12(AB+BA).
$$
Then $A^2=A\circ A=A^{\circ 2}$. Thus defined algebra is commutative but nonassociative. But one can easily verify that thus defined multiplication satisfies the  Jordan identity
$$
A^2\circ (A\circ B)=A\circ (A^2\circ B).
$$
An important property of a Jordan algebra is its \textit{power associativity}, i.e. any power of a single element $A$ is well-defined and does not depend on associations. It is not difficult to see that any element $A$ satisfies a polynomial identity (an analogue of the Cayley-Hamilton polynomial) of degree $\le n$. The rank of $A$ is just the smallest possible degree of the annihilating polynomial. One can define the trace of $A$ and the generic norm of $A$ which is just the determinant of $A$. For a general Jordan algebra, the construction is in the same spirit but more involved. The reader is referred to \cite{Schafer}, \cite{FKbook} and \cite{McCrbook} for a more detailed explanation  of the concepts discussed here and below.

Therefore  one  led to the following natural question:
\textit{What conceptually the analytical structure of Cartan isoparametric cubics (emerging from certain  PDEs) has to do with  a degree three Jordan algebra  structure}?

One possible answer is suggested in my paper \cite{Tk14}, where a one-to-one cor\-respondence between cubic homogeneous polynomial solutions  to a general eiconal equation (ge\-nera\-lizing \eqref{Muntzer0}) and degree 3 semi-simple Jordan algebras was defined. Let us briefly explain some underlying ideas of this correspondence.

 First, one associate to any  triple consisting of a vector space $V$,  a cubic form  $u$ and a nondegenerated bilinear form $\scal{}{}$ on $V$,  a commutative nonassociative algebra structure $V(u)$ in such a way that the  algebra multiplication is \textit{weakly associative}. The latter means that the identity
 \begin{equation}\label{xyz}
 \scal{xy}{z}=\scal{x}{yz}
 \end{equation}
 holds for any $x,y,z\in V$. More precisely, the multiplication $(x,y)\to xy$ is simply defined as the unique element $xy$ satisfying
 $$
 \scal{xy}{z}=u(x;y;z)
 $$
 for all $z\in V$, where $u(x;y;z)$ is the complete linearization of $u$, i.e. the unique symmetric linear form satisfying
 $$
 u(x;x;x)=6u(x).
 $$
 Thus defined multiplication is obviously commutative and satisfies \eqref{xyz}. An elementary but crucial corollary of \eqref{xyz} is that the multiplication operator $L_x$ is self-adjoint, thus, it has a real spectrum and diagonalizable.

 Recall that an idempotent of an algebra $\mathbb{A}$ is an element $c$ with the property $c^2=c$. The spectrum of $L_c$ is called the \textit{spectrum of the idempotent} $c$. To any idempotent one associates the orthogonal decomposition (called the \textit{Peirce decomposition})
 $$
 \mathbb{A}=\bigoplus_{i=1}^k\mathbb{A}_c(\lambda_i)
 $$
 of $\mathbb{A}$ into eigensubspaces $\mathbb{A}_c(\lambda_i)$, i.e. $L_c=\lambda_i$ on $\mathbb{A}_c(\lambda_i)$. If the bilinear form $\scal{}{}$ is weakly associative and positive definite, one can prove by a variational argument (see \cite{Tk18a}) that the algebra $\mathbb{A}$  necessarily contains nonzero idempotents. Indeed, any idempotent in such an algebra is (proportional to) a stationary point of the cubic form $u(x)$ on the unit sphere $\scal{x}{x}=1$ which is a compact set.

 The weak-associativity \eqref{xyz} plays a prominent role in the Jordan algebra theory, with the associating bilinear form being the generic trace form
 $$
 \tau(x,y)=\trace L_xL_y,
 $$
 where
 $$
 L_x:y\to xy
 $$
 is the (left) multiplication operator on an algebra.
 However, there is a crucial difference between the two constructions.  In the Jordan algebra theory, the trace form $\tau$ is determined by the multiplication structure and the Jordan identity
 $$
 [L_{x^2},L_{x}]=L_{x^2}L_{x}-L_{x}L_{x^2}=0
 $$
 implies that the trace form is weak associative, see \cite{McCrbook}, \cite{FKbook}. On the contrary,  the product in the algebra $V(u)$ of a cubic form $u$ is recovered from $u$ by virtue of the inner product of $V$ given by the bilinear form $\scal{}{}$ such that the latter form \textit{becomes} weak associative by its very definition. Therefore, the algebra $V(u)$ does not satisfy \textit{a priori} any identity (like Jordan algebra identity).

A correspondence  between the analytical and the algebraic sides of $V(u)$ is obtained as follows. It follows from ythe above definitions that the multiplication in $V(u)$ is completely determined  by the (polarization of) the gradient or the Hessian of $u$. More precisely,
$$
xy=D^2u(x)y.
$$
Then the Euler homogeneity theorem implies
$$
xx=x^2=D^2u(x)x=(\deg u)\nabla u(x)=2\nabla u(x),
$$
i.e. the gradient of $u$ at $x$ is $x^2$ up to a constant factor.

Thus, any PDE relation on $u$ immediately gives rise to an algebra identity on $V(u)$. This is another important ingredient of our approach. For example, Applying this construction to the Hsiang equation \eqref{Muntzer0}, and using the fact that the gradient $\nabla u(x)=x^2/2$ one obtains
\begin{equation}\label{eic}
\scal{x^2}{x^2}=36\scal{x}{x}^2
\end{equation}

With such an algebra structure in hand one is able to completely classify  cubic solutions a general eiconal equation (generalizing \eqref{Muntzer0} for an arbitrary Riemannian structure) and rely it to Jordan algebras. Then the Freudenthal-Springer-McCrimmon \cite{McCrbook} construction of an Jordan algebra from an admissible cubic form bridges these different contexts. The most nontrivial part of the latter construction is to identify and to establish that the a certain one-rank perturbation of a cubic solution  of \eqref{Muntzer0} is actually admissible. The most cubic forms \textit{are not} admissible,  but if one is lucky and such a form is found, the rest is just to follow a certain elementary algorithm \cite[p.~77]{McCrbook}.

\section{Hsiang algebras}
The above correspondence has proved very effective for the 1st order PDE, eiconal equation, and thus, for all exceptional Hsiang algebras of \textit{Cartan type}.
This opened a door for  a purely algebraic approach to tackle the general Hsiang problem.  To extend the approach on  an arbitrary Hsiang eigencubic, some further ingredients and techniques were needed. The most important part comes also from Jordan algebra theory and is known as  the Peirce decomposition. It allows to identify a finer algebra structure by virtue of multiplication rules between the so-called Peirce subspaces (the eigensubspaces of the multiplication operator by an idempotent). This method goes back to classification of associative algebras (hypercomplex number systems) in the work of the greatest nineteenth-century American mathematician Benjamin Peirce (1809-80) \cite{Peirce}. In 1930s, the method has been masterly applied  by P.~Jordan, J.~von Neumann and E.~Wigner in their seminal work on formally real Jordan algebras \cite{JordanNeumann} and developed further in a series of landmark papers of Adrian~Albert \cite{Albert42}.

Applying the above construction to the Hsiang equation \eqref{radial} yields similarly
\begin{equation}\label{L6}
\scal{x^2}{x^2}\trace L_x-\scal{x^2}{x^3}=\half{2}{3} \theta\scal{x}{x}\scal{x^2}{x}.
\end{equation}

A \textit{Hsiang algebra} is by definition any commutative algebra $\mathbb{A}$ with a weakly associative positively definite form $\scal{}{}$ and satisfying \eqref{L6}.

The correspondence between a Hsiang algebra and a solution of \eqref{radial} is given by
$$
u(x)=\frac16\scal{x}{x^2}.
$$
Thus, any Hsiang algebra produces a solution of \eqref{radial}, and conversely, any solution $u$ determines a Hsiang algebra $\mathbb{A}=V(u)$.

A Hsiang algebra is said to be trivial if its multiplication has rank one, i.e. $\dim \mathbb{A}\mathbb{A}=1$. On the level of cubic forms this means that the cubic form $u$ is essentially one-dimensional: $u(x)$ is $Cx_1^3$ in some orthogonal coordinates.

We describe very shortly some basic ideas and results following \cite[Chapter 6]{NTVbook}, \cite{Tk18e} (the full account can be found in the unpublished preprint \cite{Tk16}).

The first step is to show that for any nontrivial Hsiang algebra $\mathbb{A}$ there holds $\trace L_x=0$, i.e. the cubic form $u(x)$ is harmonic. The proof relies on the so-called minimal (extremal) idempotents, see for example \cite{Tk18b}. More precisely, one can prove that any algebra carrying a positive definite weakly associative bilinear form always has nonzero idempotents. Furthermore, for any  idempotent $c$ of the minimal possible length, the nontrivial part of its spectrum (i.e. the spectrum on the orthogonal complements to the one-dimensional eigensubspace spanned by $c$) is a subset of $(-\infty,\frac12]$. Combining this with the defining relation \eqref{L6} one is able to show the trace free property.

Now, recall that the classical Peirce approach works well for algebra identities of degree at most \textit{three} or  for at most \textit{three} different Peirce numbers. For example, coming back to isoparametric Cartan cubic algebras, the Peirce spectrum contains essentially two  elements (distinct from the trivial eigenvalue $1$) with simple fusion rules. This makes it easy to get a complete classification (see \S~3 in \cite{Tk19a}). A similar situation occurs for classical (formal real) Jordan algebras or axial algebras of Jordan type \cite{HRS15}, \cite{HSY18}. Axial algebras with $\eta=\frac12$ have a singular behaviour and require more work. Note also that the Peirce value $\frac12$ appears and plays a crucial role in the classification of nonassociative algebras associated with Hsiang exceptional eigencubics. In fact, this Peirce number is remarkable in many aspects and indicates that an algebra must satisfy to a specific algebra identity, see \cite{Tk18b}, \cite{KrTk18a},\cite{Tk18b}.

 For identities of degree higher than three, like \eqref{L6}, the situation is completely different. Some examples of such algebras are structurable algebras \cite{AlFaul}, baric algebras of degree \cite{NourigatIII}, and Majorana algebras of the Monster type \cite{Ivanov09}. In  the latter case, the algebras have four distinct Peirce numbers whose fusion rules have a natural $\mathbb{Z}_2$-grading. The full classification of all subalgebras is an important project with many potential applications to finite simple group theory, see \cite{Ivanov2018}.

But, for the Hsiang algebras the situation is more complicated because their Peirce structure is not graded. Still, it has some very nice properties. First, all idempotents in a Hsiang algebra has the same length, and, what is more important, the same Peirce spectrum. The latter property is rather extraordinary, because there are no other examples of this kind except Hsiang algebras known.  The Peirce triple
$$
(n_1,n_2,n_3)=(\dim V_c(-1),\dim V_c(-\frac12),\dim V_c(\frac12))
$$
essentially determine the structure of a Hsiang algebra.

Furthermore, given an idempotent $c\in V$, there are  two Peirce subspaces which are  \textit{subalgebras} of the ambient algebra.
It turns out that one can deform the original multiplication in these subalgebras such that some \textit{hidden} algebra structures become visible: a Clifford algebra structure  on $V_c(1)\oplus V_c(-1)$ (constructed ad hoc) and a formally real Jordan algebra structure $\Lambda_c$ on $V_c(1)\oplus V_c(-\frac12)$ (obtained by using the Freudenthal-Springer-McCrimmon construction). Manipulating with this structures in hands, one can establish the basic classification including \ref{item1}--\ref{item6} by pure algebraic means.

To absorb exceptional algebras, one establishes the following result:  the Jordan algebra structure $\Lambda_c$ on $V_c(1)\oplus V_c(-\frac12)$ is that it is simple if and only if the ambient Hsiang algebra is exceptional. The simple formally real Jordan algebras are well-known since the celebrated classification of P.~Jordan, J.~von Neumann and E.~Wigner  Jordan algebras \cite{JordanNeumann}. Since one also has an additional obstruction coming form the existence of Clifford algebra structure  on $V_c(1)\oplus V_c(-1)$, one is able to  deduce the finiteness of  admissible Peirce triples for exceptional algebras. Those are displayed in Table~\ref{tabs} below.

\begin{footnotesize}
\def\MM{6mm}
\begin{table}[ht]
\renewcommand\arraystretch{1.5}
\noindent
\begin{flushright}
\begin{tabular}{p{10mm}|p{6.5pt}|p{6.5pt}|p{6.5pt}|p{6.5pt}|
p{6.5pt}|p{6.5pt}|p{6.5pt}|p{6.5pt}|p{6.5pt}|p{6.5pt}|p{6.5pt}|p{6.5pt}|p{6.5pt}|p{6.5pt}|p{6.5pt}|p{6.5pt}|p{6.5pt}|p{6.5pt}|p{6.5pt}|
p{6.5pt}|p{6.5pt}|p{6.5pt}|p{6.5pt}}
$\dim V$  &  $5$ & $8$  & $14$  & $26$  &$9$ & $12$ & $\textcolor{black}{15}$ & $21$ & $15$ & $18$ & $\textcolor{black}{21}$ & $\textcolor{black}{24}$ &  $\textcolor{black}{30}$&  $\textcolor{black}{42}$& $27$ & $30$  & $\textcolor{black}{33}$  & $\textcolor{black}{36}$  & $\textcolor{black}{51}$ & $54$ & $\textcolor{black}{57}$ & $\textcolor{black}{60}$ & $\textcolor{black}{72}$      \\\hline
$n_1$ &$2$ & $3$  & $5$  & $9$  & $0$ & $1$  & $2$  & $4$  & $0$  & $1$  &  $2$ & $3$ &  $5$ & $9$  & $0$  & $1$   &  $2$  & $3$   & $0$  & $1$  & $2$ & $3$ & $7$ \\\hline
$n_2$ &  $0$ & $0$  & $0$  & $0$  &$5$ & $5$  & $5$  & $5$  & $8$  & $8$  &  $8$ &  $8$ &  $8$ & $8$  &$14$  & $14$  & $14$  & $14$  & $26$ & $26$ & $26$ & $26$ & $26$\\\hline
$n_3$ & $2$ & $4$  & $8$  & $16$  &$3$ & $5$  & $7$  & $11\,$ & $6$  & $8$  & $10\,$ & $12\,$ & $16\,\,\,$ & $24\,$ &$12\,$  & $14\,$  & $16\,\,$  & $18\,$  & $24\,$ & $26\,$ & $28\,$ & $30\,\,$ & $38\,\,$      \\\hline
\end{tabular}
\end{flushright}
\bigskip
\caption{The a priori admissible dimensions $\dim V$ and the Peirce triples $(n_1,n_2,n_3)=(\dim V_c(-1),\dim V_c(-\frac12),\dim V_c(\frac12))$ of exceptional algebras}\label{tabs}
\end{table}
\end{footnotesize}

Thus, the proposed program also yield some further results on Hsiang eigencubics, but an ultimate classification required a deeper insight into the structure of exceptional eigencubics. The most difficult part of the program was (and is) to eliminate `false' exceptional Hsiang algebras from the list predicted by \ref{item3} and displayed in Table~\ref{tabs}. To this aim, I developed in 2014 a next important ingredient, the so-called tetrad decomposition. This technique seems to me especially  important for some further applications even beyond the Hsiang problem. In short, the idea is as follows. We already seen above that the Peirce subspace $V_c(1)\oplus V_c(-\frac12)$ can be deformed into a Jordan algebra $\Lambda_c$. The algebra multiplication on $\Lambda_c$ is different from (but isotopic to) the original multiplication on $V$; more precisely it is a rank one deformation of the latter. The Jordan structure on $\Lambda_c$ makes visible some hidden structures in $V$. For example, the idempotent $c$ is the unit in $\Lambda_c$  and any nontrivial (distinct from zero and unit) idempotent in $\Lambda_c$ is an absolute nilpotent (square zero elements) in $V$,  and vice versa. On the Jordan algebra level it is well-known that for any primitive idempotent (i.e. one which cannot be written as a sum of non-zero  idempotents) in $\Lambda_c$ gives rise to a \textit{Jordan frame}. The latter is just a partition of the Jordan algebra unit into orthogonal primitive idempotents (played an important tool in the spectral theory of Jordan algebras by P.~Jordan, J.~von Neumann and E.~Wigner). This establishes a natural connection  these two structures.

By using the tetrad decomposition, one can show that some Peirce dimensions in the above table are not realizable. This eliminates the Perice triples $(n_1,n_2,n_3)$ with values $(2,5,7)$, $(2,8,10)$, $(2,14,16)$, $(3,14,18)$, $(0,26,24)$, $(2,26,28)$, $(3,26,30)$ and $(7,26,38)$. All triples with $n_2=0$, $n_1=0$ and $n_2\ne26$, $n_1=1$ and $(n_1,n_2,n_3)=(4,5,11)$ are realizable, see the description below:

\begin{itemize}
\item
If ${n_2 =0}$ then the corresponding Hsiang algebras are the contractions of the algebra of the cubic form $u=\frac{1}{6}\scal{z}{z^2}$ on the trace free subspace  of the Jordan algebra $\mathscr{H}_3({\Bbb K}_d)$ of $3\times 3$ Hermitian matrices over a real division algebra ${\Bbb K}_d$ of dimension $d=1,2,4,8$.

\medskip
\item
If ${n_1=0}$ then $n_2\in \{5,8,14\}$ and the corresponding Hsiang algebras are the algebras of cubic forms $\frac{1}{12}\scal{z^2}{3\bar z-z}$, where $z\to \bar z$ is the natural involution on $V= \mathscr{H}_3({\Bbb K}_d)$, $d=2,4,8$.

\medskip
\item
If ${n_1=1}$ then $n_2\in \{5,8,14,26\}$ and the corresponding Hsiang algebras are the algebras of cubic forms  $u(z)=\re \scal{z}{z^2}$ on the complexification $\mathscr{H}_3({\Bbb K}_d)\otimes \mathbb{C}$, $d=1,2,4,8$.

\medskip
\item
If ${(n_1,n_2)=(4,5)}$ then $V$ is the algebra of cubic form $u=\frac{1}{6}\scal{z}{z^2}$ on the contraction of the Albert algebra on the purely imaginary subspace: $\mathscr{H}_3({\Bbb K}_8)\ominus \mathscr{H}_3({\Bbb K}_1)$.

\end{itemize}

The three remained triples with  $n_2=8$ are still an open problem.

\section{Conclusion}

When I returned from military service in 1988, Miklyukov was completely passionate about the string theory and its  connections to zero mean curvature surfaces in Lorenz spaces which has been his principal and very fruitful  direction during the 1990s. Together with his students (primarily the twin brothers Vladimir Klyachin and Alexei Klyachin) he published a series of very deep existence and regularity results on maximal surfaces.

We together came back to our project on harmonic functions on Riemannian manifolds and external structure of minimal surfaces later, in the mid of 1990s. We were able to relax the entire graph condition in the Bernstein result for slowly growing multiplicities by using the fundamental frequency technique. This became our third and the last join paper \cite{MT96}.

The fundamental frequency technique is very powerful and not yet fully used method with many potential applications. For example, in 2005 I  developed some of our ideas in \cite{MT87} and \cite{MT96} to apply them to the Meeks problem on disjoint minimal graphs. The main result of \cite{Tk09a} states that there may be at most 3 minimal graphs supported on disjoint domains of $\R{2}$ (Meeks and Rosenberg conjectured in \cite{Meeks} that the optimal number must be 2, which remains  an open problem). Both the two- and higher-dimensional cases have natural connections to the Bernstein property. A better understanding of this point in the 2D and higher-dimensions is an ongoing project with Luciano Mari, SNS, Pisa.

The ultimate classification of nonassocaitive algebras of cubic minimal cones is still an incomplete project but a complete picture seems clear. A interesting and deep direction here is a better understanding of general metrisable algebras and their Peirce structure in different differential geometric and group theoretic contexts \cite{KrTk18a}, \cite{Fox19a}, \cite{HRS15b}, \cite{Ivanov2018}.

Coming back to the central question about  a possible connection between the existence of Milnor's spheres and the breakdown of the Bernstein property, we have to recognize that we still don't understand the essence of the matter. On the other hand, the algebraic structures appeared in the study of minimal cones have a very close relation to the algebraic part of the Milnor construction. Perhaps, the answer could be found from a better understanding of the unusual (H\"older continuous) viscosity solutions constructed very recently by Nikolai Nadirashvili and Serge Vl{\u{a}}du{\c{t}} \cite{NV07},  \cite{NV08}, \cite{NV11a}, \cite{NV17}, and also in \cite{NTV}, \cite{NTVbook} by using certain nonassociative algebra structures. Remarkably, the zero varieties of these solutions are minimal cones but  no conceptional explanation for this fact is known so far.

\bibliographystyle{plain}
\def\cprime{$'$}

\end{document}